\documentclass[preprint,12pt]{elsarticle}

\usepackage{amssymb}
\usepackage{amscd,amsmath,amsfonts,graphicx}
\usepackage{amsthm}

\begin{document}
\baselineskip=22pt \centerline{\large \bf Characterizations of the
Logistic and Related Distributions} \vspace{1cm}
\centerline{Chin-Yuan Hu and Gwo Dong Lin} \centerline{National
Changhua University of Education and Academia Sinica} \vspace{1cm}
\noindent {\bf Abstract.} It is known that few characterization
results of the logistic distribution were available before, although
it is similar in shape to the normal one whose characteristic
properties have been well investigated. Fortunately, in the last
decade, several authors have made great progress in this topic. Some
interesting characterization results of the logistic distribution
have been developed recently. In this paper, we further provide some
new results by the distributional equalities in terms of order
statistics of the underlying distribution and the random exponential
shifts. The characterization of the closely related Pareto type II
distribution is also investigated.
\vspace{0.1cm}\\
\hrule
\bigskip
\noindent AMS subject classifications: Primary 62E10, 62G30, 60E10.\\
\noindent {\bf Key words and phrases:} Characterization, order
statistics, stochastic order, logistic distribution,
exponential distribution, Pareto type II distribution.\\
{\bf Short title: The logistic and related distributions}\\
{\bf Postal addresses:} Chin-Yuan Hu, Department of Business
Education, National Changhua University of Education, Changhua
50058, Taiwan. (E-mail: buhuua@gmail.com)\\
Gwo Dong Lin, Institute of Statistical Science, Academia Sinica,
Taipei 11529, Taiwan. (E-mail:
 gdlin@stat.sinica.edu.tw)
\newpage

\noindent{\bf 1. Introduction}
\newcommand{\bin}[2]{
   \left(
     \begin{array}{@{}c@{}}
         #1  \\  #2
     \end{array}
   \right)          }

The logistic distribution is similar to a normal distribution in
shape (Mudholkar and George 1978) and has an explicit closed form,
so it has some advantages
 in practical applications.  As remarked by Kotz (1974), few characterizations of  the logistic distribution
 were available before, but recently, some interesting
results  have been developed. In this paper, we will further provide
some more new results by properties of order statistics.

 \indent We first introduce some notations. Let $X$ obey the
distribution $F$, denoted by $X\sim F.$ Let $\{X_j\}_{j=1}^n$ be a
random sample of size $n$ from distribution $F$ and denote  the
corresponding order statistics by $X_{1,n}\leq X_{2,n}\leq \cdots
\leq X_{n,n}$. The distribution function of $X_{k,n}$ is denoted by
$F_{k,n}$. It is known that $F_{k,n}$ is the composition of
$B_{k,n-k+1}$ and $F$ (see, e.g., Hwang and Lin 1984), where
$B_{\alpha,\beta}$ is the beta distribution with parameters
$\alpha,\beta>0$, namely,
\begin{eqnarray}F_{k,n}(x)&=&B_{k,n-k+1}(F(x))=k\bin{n}{k}\int_0^{F(x)}t^{k-1}(1-t)^{n-k}dt,\ \ x\in
{\mathbb R}\equiv(-\infty,\infty),\\
B_{\alpha,\beta}(u)&=&\frac{\Gamma(\alpha+\beta)}{\Gamma(\alpha)\Gamma(\beta)}\int_0^ut^{\alpha-1}(1-t)^{\beta-1}dt,\ \
u\in[0,1].
\end{eqnarray}

On the other hand, for $Y\sim G,$ we say that $X$ is less than or
equal to $Y$ in the usual stochastic order, denoted by $X\le_{st}
Y,$ if $\overline{F}\le \overline{G},$ where
$\overline{F}(x)=1-F(x)=\Pr(X>x).$

Let us start with an interesting  simple example. Clearly, for
general distribution $F$, we have $X_{1,2}\le_{st} X$ because
${F}_{1,2}=1-\overline{F}^2\ge {F}$ by (1). One possibility to
adjust this "inequality"  is to choose a nonnegative random variable
$Z$, independent of $X$ and $X_j$'s, such that
\begin{eqnarray}X\stackrel{\rm d}{=}X_{1,2}+Z
\end{eqnarray}
or
\begin{eqnarray}X-Z\stackrel{\rm d}{=}X_{1,2},
\end{eqnarray}
where $\stackrel{\rm d}{=}$ means equality in distribution. One
might think that the solutions of the distributional equations (3)
and (4) are the same, but this is not true in general, because the
characteristic function of $Z$ is not equal to the reciprocal of
that of $-Z,$ namely, $E[e^{itZ}]\ne (E[e^{-itZ}])^{-1},\ t\in
{\mathbb R},$ in general. For example, if $Z$ has the standard
exponential distribution ${\cal E}$, then the solution of (3) is a
logistic distribution $F(x)=1/[1+e^{-(x-\mu)}],\ x\in{\mathbb R},$
where $\mu\in {\mathbb R}$ is a constant, while the solution of (4)
is a negative (or reversed) exponential distribution
$F(x)=e^{(x-\mu)/2},\ x\le \mu,$ where $\mu\in {\mathbb R}$ is a
constant (see, e.g., George and Mudholkar 1982, Lin and Hu 2008, and
Ahsanullah et al.\,2011, and note that the smoothness conditions on
$F$ therein are redundant due to Lemmas 1--3 below).

Throughout the paper, let $U$ and $\xi$ obey the uniform
distribution ${\cal U}$ on $[0,1]$ and  the standard exponential
distribution ${\cal E},$ respectively. Moreover, let
$\{U_j\}_{j=1}^n$ and $\{U_j^{\prime}\}_{j=1}^n$ be two random
samples of size $n$ from  ${\cal U}$, and let $\{\xi_j\}_{j=1}^n$
and $\{\xi_j^{\prime}\}_{j=1}^n$ be two random samples of size $n$
from ${\cal E}.$
 All the above random variables $X$, $U$, $\xi$,
$X_j$, $U_j$, $U_j^{\prime}$, $\xi_j$ and $\xi_j^{\prime},
j=1,2,\ldots,n,$ are assumed to be independent from now on.

Mimicking the above characterization approaches (3) and (4), several
authors have considered the general stochastic inequality
$X_{k,n}\le_{st} X_{k+1,n}$ and solved the distributional equations
(a) $X_{k+1,n}\stackrel{\rm d}{=}X_{k,n}+a\xi$ and (b)
$X_{k,n}\stackrel{\rm d}{=} X_{k+1,n}-b\xi$, or, more generally, (c)
$X_{k,n}+a\xi_1\stackrel{\rm d}{=} X_{k+1,n}-b\xi_2$, where $a$ and
$b$ are nonnegative constants. In particular, the equality
$$X_{k,n}+\frac{1}{n-k}\xi_1\stackrel{\rm d}{=} X_{k+1,n}-\frac{1}{k}\xi_2 $$
also characterizes the logistic distribution.  (See AlZaid and
Ahsanulla 2003, Weso$\hbox{\l}$owski and Ahsanulla 2004, and
Ahsanulla et al.\,2010, 2012 for equations arising from
$X_{k,n}\le_{st} X_{k+r,n}$ with $1\le r\le n-k$.) Besides, the
distributional equations arising from (a) $X\le_{st} X_{n,n}$, (b)
$X_{1,n}\le_{st} X_{n,n}$ and (c) $X_{m,m}\le_{st} X_{n,n}$, where
$m<n$,  were investigated by Zykov and Nevzorov (2011), Ananjevskii
and Nevzorov (2016) as well as Berred and Nevzorov (2013),
respectively.
\medskip\\
\indent
In this paper we will solve the distributional equations arising
from stochastic inequalities: (i) $X_{k,n}\le_{st} X_{k,n-1}$, (ii)
$X_{k,n-1}\le_{st} X_{k+1,n}$, (iii) $X_{k,n}\le_{st} X_{k,k}$, (iv)
$X_{1,k}\le_{st} X_{n-k+1,n}$,  (v) $X_{k,n}\le_{st} X_{k,n-m}$, and
(vi) $X_{m-k,n-k}\le_{st} X_{m,n}.$
\medskip\\
\indent To do this, some useful lemmas are given in the next
section. The main characterization results are stated and proved in
Sections 3 and 4. For simplicity, we first deal with the closely
related Pareto type II distribution in Section 3, and then the
logistic distribution in Section 4. Finally, we pose an open problem
in Section 5.
\bigskip\\
 \noindent{\bf 2. Lemmas}

 We need some lemmas in the sequel. Lemma 1(i) was given without proof in Lukacs (1970, p.\,38), but
 has been ignored in the literature.
 We provide a proof here for completeness.\\
\noindent{\bf Lemma 1.} (i) Let $Y$ and $Z$ be two independent
random variables. If $Y$ has an absolutely continuous distribution,
then so does $Y+Z,$ regardless of the distribution of $Z.$\\
(ii) If, in addition to the assumptions in (i), both $Y$ and $Z$ are
positive random variables, then the product $YZ$ has an absolutely
continuous distribution.\\
{\bf Proof.} Let $F$, $G$ and $H$ be the distributions of $Y+Z$, $Y$
and $Z$, respectively. Then
$$F(x)=\int_{-\infty}^{\infty}G(x-z)dH(z),\ \ x\in{\mathbb R}.$$
Since $G$ is absolutely continuous, we have that for each
$\varepsilon>0$, there exists a $\delta>0$ such that
$\sum_{j=1}^n[G(b_j)-G(a_{j})]<\varepsilon$\, if\, $\sum_{j=1}^n
(b_j-a_j)<\delta,$ where $a_j<b_j\le a_{j+1}<b_{j+1},\
j=1,2,\ldots,n-1.$ This in turn implies that for the above $\{(a_j,
b_j)\}_{j=1}^n$,
$$\sum_{j=1}^n[F(b_j)-F(a_{j})]
=\int_{-\infty}^{\infty}\sum_{j=1}^n[G(b_j-z)-G(a_j-z)]dH(z)<\int_{-\infty}^{\infty}\varepsilon\ dH(z)=\varepsilon.
$$
Hence, part (i) is proved. To prove part (ii), we recall first that
both the logarithmic and exponential functions are absolutely
continuous, and that the composition preserves the property of
absolute continuity. Then consider $\log(YZ)=\log Y+\log Z$ and use
part (i) to conclude that $\log(YZ)$ has an absolutely continuous
distribution, and hence so does $YZ.$ This completes the proof.

 It is known that the inverse function of an absolutely continuous function
 with positive derivative {\it almost everywhere} is not
 necessarily absolutely continuous. However, we have the following useful result.

\noindent{\bf Lemma 2.} Let $F$ be an absolutely continuous
distribution on $[0,1]$ and $F^{\prime}(x)=f(x)>0$ on $(0,1).$ Then
the inverse function of $F$ is itself
an absolutely continuous distribution.\\
{\bf Proof.} By the assumptions, $F$ is a strictly increasing and
continuous function from $[0,1]$ to $[0,1]$, so is its inverse
function $F^{-1}.$ This implies that $F^{-1}$ is a continuous
distribution on $[0,1].$ Moreover,
$\frac{d}{dt}F^{-1}(t)=1/f(F^{-1}(t))$ is a positive measurable
function on $(0,1)$ (see, e.g., Shorack and Wellner 1986, pp.\,8-9).
By changing variables $x=F^{-1}(t),$ we have
$\int_0^1[\frac{d}{dt}F^{-1}(t)]dt=\int_0^1[1/f(x)]\cdot f(x)dx=1.$
Therefore, the distribution function $F^{-1}$ has no singular part
and is absolutely continuous. The proof is complete.

 \noindent{\bf Lemma 3.} If the distribution $F_{k,n}$
of order statistic $X_{k,n}$ is absolutely continuous, then so is
the underlying distribution $F$ of $X.$\\
{\bf Proof.} Recall (see (1)) that $F_{k,n}(x)=B_{k, n-k+1}(F(x)),\
x\in{\mathbb R},$ where $B_{k, n-k+1}$, defined in (2), is the beta
distribution $B_{\alpha,\beta}$ with parameters $\alpha=k$ and
$\beta=n-k+1$, and has a positive continuous density function on
$(0,1).$ Therefore, $F=B_{k, n-k+1}^{-1}\circ F_{k,n}$ is absolutely
continuous by Lemma 2. The proof is complete.

\noindent{\bf Lemma 4.} Let $1\le k<n.$ Then we have the following identities:\\
(i) $F_{k,n}-F_{k+1,n}=\bin{n}{k}F^k{\overline{F}}^{n-k},$\ \
(ii) $F_{k,n}-F_{k,n-1}=\bin{n-1}{k-1}F^k{\overline{F}}^{n-k},$ and\\
(iii) $F_{k,n-1}-F_{k+1,n}=\bin{n-1}{k}F^k{\overline{F}}^{n-k}.$\\
{\bf Proof.} For parts (i) and (ii), see, e.g.,  David and Shu
(1978) as well as David and Nagaraja (2003, p.\,23), while part
(iii) follows from parts (i) and (ii) and the identity:
$$\bin{n}{k}=\bin{n-1}{k-1}+\bin{n-1}{k}.$$

\indent Denote  the left and the right extremities of
$F$ by $\ell_F$ and $r_F$, respectively. It is known that if the
absolutely continuous $F$ satisfies the functional equation
$F^{\prime}(x)=F(x)(1-F(x)),\ x\in(\ell_F,r_F),$ then $F$ is a
logistic distribution. In the next lemma we extend this result.

 \noindent{\bf Lemma 5.} Let $r,\theta>0,$ $a\in[0,1],$ and let $F$ be an
absolutely continuous distribution function satisfying\ \
$x^aF^{\prime}(x)=\theta F(x)(1-F^r(x))$ for $x\in(\ell_F,r_F).$\\
(i) If $a=1,$  then $\ell_F=0, r_F=+\infty$ and
$$F(x)=\left(\frac{\lambda x^{r\theta}}{1+\lambda
x^{r\theta}}\right)^{1/r},\ \ 0\le x<\infty,$$ where $\lambda$ is a
positive constant.\\
(ii) If  $a\in[0,1),$ then $\ell_F=-\infty, r_F=+\infty$ and
$$F(x)=\left(\frac{\lambda \exp({\frac{r\theta}{1-a}x^{1-a}})}{1+\lambda \exp({\frac{r\theta}{1-a}x^{1-a}})}\right)^{1/r},\ \ -\infty< x<\infty,$$
 where $\lambda$ is a positive constant.\\
{\bf Proof.} Define the increasing function $G(x)=(1-F^r(x))^{-1}-1$
from $(\ell_F,r_F)$ onto $(0,\infty).$ Then
$G^{\prime}(x)=rF^{r-1}(x)F^{\prime}(x)(1-F^r(x))^{-2},$ and hence
$${x^{a}}G^{\prime}(x)={r\theta}G(x),\ x\in(\ell_F,r_F).$$
(a) If $a=1,$  solving the above equation leads to
 $G(x)=\lambda
x^{r\theta},\ x\in(\ell_F,r_F),$ for some constant $\lambda>0.$ On
the other hand, we have, by the definition of $G$, that
$$F(x)=\left(\frac{G(x)}{1+G(x)}\right)^{1/r}=\left(\frac{\lambda x^{r\theta}}{1+\lambda
x^{r\theta}}\right)^{1/r},\ \ x\in(\ell_F,r_F),$$ and hence,
$\ell_F=0$ and $r_F=+\infty,$ because $F$ is a distribution
function. This proves part (i).\\
(b) If $a\in [0,1),$  we have instead $G(x)=\lambda
\exp({\frac{r\theta}{1-a}x^{1-a}}),\ x\in(\ell_F,r_F),$ for some
constant $\lambda>0.$ The required result then follows from both the
definition of the function $G$ and  the fact that $F$ is a
distribution function. The proof is complete.

 Some equalities (in
distribution) of the next lemma are essentially due to Nevzorov
(2001, Lecture 3), but we provide here an alternative and possibly
simpler proof.

\noindent{\bf Lemma 6.} Let $\xi_{k,n}$ ($U_{k,n}$, resp.) be the
$k$-th smallest order statistic of a random sample of size $n$ from
the standard exponential distribution ${\cal E}$ (the uniform
distribution ${\cal U}$, resp.).
Then the following statements are true.\\
(i) The Laplace transform of $\xi_{k,n}$ is
$L_{\xi_{k,n}}(s)=\frac{n-k+1}{n-k+1+s}\cdot
\frac{n-k+2}{n-k+2+s}\cdots\frac{n}{n+s},\ \ s\ge 0.$\\
(ii) $\xi_{k,n}\stackrel{\rm d}{=}\sum_{j=1}^k\frac{1}{n-j+1}\xi_{k-j+1},$ where $1\le k\le n.$\\
 (iii) $\xi_{m,n}\stackrel{\rm
d}{=}\xi_{m-k,n-k}+\xi_{k,n}^{\prime},$ where $1\le k<m\le n.$\\
(iv) $U_{n-k+1,n}\stackrel{\rm d}{=}\prod _{j=1}^kU_j^{1/(n-j+1)},$  where $1\le k\le n.$\\
(v) $U_{k,n}\stackrel{\rm d}{=}U_{k,m-1}\cdot U_{m,n}^{\prime},$
where $1\le k<m\le n.$\\
{\bf Proof.} Recall that the distribution of $\xi_{k,n}$ is
$F_{\xi_{k,n}}(x)=k\bin{n}{k}\int_0^{{\cal
E}(x)}t^{k-1}(1-t)^{n-k}dt,\ \ x\ge 0,$ where ${\cal
E}(x)=F_{\xi}(x)=1-e^{-x},\ x\ge 0.$ Then the Laplace transform of
$\xi_{k,n}$ is
$$L_{\xi_{k,n}}(s)=E[e^{-s\xi_{k,n}}]=\int_0^{\infty}e^{-sx}dF_{\xi_{k,n}}(x)=k\bin{n}{k}\int_0^{\infty}e^{-(n-k+1+s)x}(1-e^{-x})^{k-1}dx,\ \ s\ge 0.$$
By integration by parts, it follows from the above integral that
\begin{eqnarray*}
L_{\xi_{k,n}}(s)&=&k\bin{n}{k}\frac{k-1}{n-k+1+s}\int_0^{\infty}e^{-(n-k+2+s)x}(1-e^{-x})^{k-2}dx=\cdots\\
&=&k\bin{n}{k}\frac{k-1}{n-k+1+s}\cdot \frac{k-2}{n-k+2+s}\cdots\frac{1}{n-1+s}\cdot\int_0^{\infty}e^{-(n+s)x}dx\\
&=&\frac{n-k+1}{n-k+1+s}\cdot
\frac{n-k+2}{n-k+2+s}\cdots\frac{n-1}{n-1+s}\cdot\frac{n}{n+s},\ \ \ s\ge 0.
\end{eqnarray*}
This proves part (i), which in turn implies parts (ii) and (iii) by
the fact that $E[e^{-s(\xi/k)}]=k/(k+s),\ s\ge 0.$ Part (iv) follows
from part (ii) because  $U_j\stackrel{\rm d}{=}\exp(-\xi_{k-j+1})$
and the order statistic $U_{n-k+1,n}\stackrel{\rm
d}{=}\exp(-\xi_{k,n}).$ To prove part (v), we have
$U_{n-m+1,n}\stackrel{\rm d}{=}U_{n-m+1,n-k}\cdot
U_{n-k+1,n}^{\prime}$ by using part (iii), and then reset $k=n-m+1.$
The proof is complete.

\noindent{\bf Lemma 7.} Let $Y$ obey the Pareto type II (or
log-logistic) distribution $G(y)={y}/(1+y),\ y\ge 0.$ Let
$\{Y_j\}_{j=1}^n$, independent of $U$ and $\{U_j\}_{j=1}^n$,  be a
random sample of size $n$ from $G$. Then we have the following
equalities in distribution:\\
(i) $1/Y\stackrel{\rm d}{=}Y$ and in general,
$1/Y_{k,n}\stackrel{\rm
d}{=}Y_{n-k+1,n},$  where $1\le k\le n.$\\
(ii) $Y_{k,n-1}\stackrel{\rm
d}{=}Y_{k,n}/U^{1/(n-k)},$ where $1\le k\le n-1.$\\
(iii) $Y_{k,n-m}\stackrel{\rm
d}{=}Y_{k,n}/U_{n-k-m+1,n-k},$ where $1\le k\le n-m.$\\
(iv) $Y_{k,n-1}\stackrel{\rm
d}{=}Y_{k+1,n}\cdot U^{1/k},$ where $1\le k\le n-1.$\\
(v) $Y_{m-k,n-k}\stackrel{\rm d}{=}Y_{m,n}\cdot U_{m-k,m-1},$
where $2\le k+1\le m\le n.$\\
{\bf Proof.}  It is easy to check part (i). To prove the remaining
parts, recall that the distribution of $Y_{k,n}$ is
$G_{k,n}(y)=k\bin{n}{k}\int_0^{G(y)}t^{k-1}(1-t)^{n-k}dt,\ y\ge 0.$
Then we have $H(y)\equiv \Pr(Y_{k,n}/U^{1/(n-k)}\le
y)=\int_0^1G_{k,n}(yu^{1/(n-k)})du.$ By changing variables,
\begin{eqnarray*}H(y)=\int_0^1G_{k,n}(yt)dt^{n-k}=(n-k)\int_0^1G_{k,n}(yt)t^{n-k-1}dt,\ \ y\ge 0.
\end{eqnarray*}
Therefore, $G_{k,n-1}=H$ iff, by differentiation,
$$\left(\frac{1}{1+y}\right)^n=\int_0^1n\left(\frac{1}{1+yt}\right)^{n+1}t^{n-1}dt,\ \ y\ge 0,$$
which is, however, a special case of the well-known identity
$$\left(\frac{1}{1+y}\right)^{\beta_1}=\frac{\Gamma(\beta_1+\beta_2)}{\Gamma(\beta_1)\Gamma(\beta_2)}
\int_0^1\left(\frac{1}{1+yt}\right)^{\beta_1+\beta_2}t^{\beta_1-1}(1-t)^{\beta_2-1}dt,\
\ y\ge 0$$ (see Gradshteyn and Ryzhik 2014, p.\,314). This proves
part (ii). Part (iii) follows from part (ii) by iteration and Lemma
6(iv), while part (iv) follows from either parts (i) and (ii)
(letting $k_*=n-k$) or Lemma 8(iii) below because $Y\stackrel{\rm
d}{=}\exp(X)$ and $U\stackrel{\rm d}{=}\exp(-\xi).$ Finally,  we
prove part (v) by using part (iv), iteration and Lemma 6(iv) again.
The proof is complete.
\\
 \noindent{\bf
Lemma 8.} Let $X$ obey the standard logistic distribution
$F(x)=1/[1+\exp(-x)],\ x\in{\mathbb R}.$ Then we have the following
equalities in distribution:\\
(i) $X_{k,n-1}\stackrel{\rm d}{=}X_{k,n}+ \frac{1}{n-k}\xi,$ where
$1\le
k\le n-1.$\\
(ii) $X_{k,n-m}\stackrel{\rm d}{=}X_{k,n}+ \xi_{m,n-k},$ where $1\le
k\le n-m.$\\
 (iii) $X_{k,n-1}\stackrel{\rm
d}{=}X_{k+1,n}-\frac{1}{k}\xi,$ where $1\le k\le n-1.$\\
(iv) $X_{m-k,n-k}\stackrel{\rm d}{=}X_{m,n}-\xi_{k,m-1},$ where
$2\le k+1\le m\le n.$\\
{\bf Proof.} The results follow from Lemma 7  by noting that (a)
$X\stackrel{\rm d}{=}\log Y,$ (b) $\xi \stackrel{\rm d}{=} -\log U$
and (c) $-\log U_{k,n}\stackrel{\rm d}{=}\xi_{n-k+1,n}$  for all
$1\le k\le n.$  The proof is complete.

For the proof of the next lemma, see Lin and Hu (2008, Lemma 5).

\noindent{\bf Lemma 9.} Let $f$ and $g$ be two functions real
analytic and strictly monotone in $[0,\infty)$. Assume that for each
$n\geq 1$, the $n$-th derivatives $f^{(n)}$ and $g^{(n)}$ are
strictly monotone in some interval $[0,\delta_n)$. Let
$\{x_n\}_{n=1}^{\infty}$ be a sequence of positive real numbers
converging to zero. If $f(x_n)=g(x_n),
n=1,2,\ldots,$ then $f=g$.\medskip\\
\noindent{\bf 3. Characterizations of the Pareto type II
distribution}

We start with the Pareto type II distribution which is easier to
handle, and recall  that the uniform order statistic $U_{k,n}\sim
B_{k,n-k+1}.$

 \noindent{\bf Theorem 1.} Let $Y\sim G$ be a positive random variable and let $1\le k \le n-1$ be fixed integers.
 Let $Y_1,Y_2,\ldots,Y_n$ be $n$ independent copies of $Y,$ and let
$U$, independent of $\{Y_i\}_{i=1}^n$, be a random variable with
uniform distribution on $[0,1].$ Then the distributional equality
\begin{eqnarray}
Y_{k,n-1}\stackrel{\rm d}{=}Y_{k,n}/U^{1/(n-k)}
\end{eqnarray}
holds  iff $G$ is a Pareto distribution
$G(y)=\lambda y/(1+\lambda y),\ y\ge 0,$ where $\lambda>0$ is a constant.\\
{\bf Proof.} The sufficiency part follows from Lemma 7(ii) because
$(\lambda Y)_{\ell, m}\stackrel{\rm d}{=}\lambda Y_{\ell,m}$ for
$\lambda>0$ and for all $1\le\ell\le m.$ To prove the necessity
part, we note, by Lemma 1(ii), that the distribution $G_{k,n-1}$ of
$Y_{k,n-1}$ is absolutely continuous, and so is $G$ by Lemmas 2 and
3. Rewrite (5) as
$$G_{k,n-1}(y)=\int_0^{1}G_{k,n}(yu^{1/(n-k)})du,\ \ y\ge 0.$$
By changing variables $t=yu^{1/(n-k)},$ we have
$$G_{k,n-1}(y)=(n-k)y^{-(n-k)}\int_0^{y}G_{k,n}(t)t^{n-k-1}dt,\ \ y>0.$$
Taking differentiation leads to
\begin{eqnarray}
yG_{k,n-1}^{\prime}(y)=(n-k)[G_{k,n}(y)-G_{k,n-1}(y)],\ \ y>0.
\end{eqnarray}
With the help of (1) and Lemma 4(ii), (6) is equivalent to
$$yG^{\prime}(y)=G(y)[1-G(y)],\ \ y\in(\ell_G,r_G).$$
Finally, Lemma 5(i) with $r=\theta=1$ completes the proof.

\noindent{\bf Corollary 1.}  Under the same assumptions of Theorem
1, the distributional  equality $Y_{k,n-1}\stackrel{\rm
d}{=}Y_{k,n}/U^{1/\alpha}$ holds for some $\alpha>0$, iff $G$ is a
Pareto distribution with $\overline{G}(y)=1/[1+\lambda
y^{\alpha/(n-k)}],\ y\ge 0,$ where $\lambda$ is a positive constant.\\
{\bf Proof.} Consider $Y_i^{\prime}=Y_i^{\alpha/(n-k)}$ for
$i=1,2,\ldots.$ Then $Y^{\prime}_{k,n}=Y_{k,n}^{\alpha/(n-k)},$ and
apply Theorem 1 to  the case: $Y_{k,n-1}^{\prime}\stackrel{\rm
d}{=}Y_{k,n}^{\prime}/U^{1/(n-k)}.$

\noindent{\bf Corollary  2.} Under the same assumptions of Theorem
1, the distributional  equality
\begin{eqnarray*}
Y_{k,n-1}\stackrel{\rm d}{=}Y_{k+1,n}\cdot U^{1/k}
\end{eqnarray*}
holds  iff $G$ is a Pareto distribution $G(y)=\lambda y/(1+\lambda
y),\ y\ge 0,$  where
$\lambda>0$ is a constant.\\
{\bf Proof.} The sufficiency part is a consequence of Lemma 7(iv).
 To prove the necessity
part, denote $Y^*=1/Y.$ Then $Y^*_{\ell, m}=(1/Y)_{\ell,
m}\stackrel{\rm d}{=}1/ Y_{m-\ell+1,m}$ for all $1\le\ell\le m.$ By
assumptions, we have the equality $1/Y_{k,n-1}\stackrel{\rm
d}{=}1/Y_{k+1,n}\cdot 1/U^{1/k},$ or, equivalently,
$Y^*_{n-k,n-1}\stackrel{\rm d}{=}Y^*_{n-k,n}\cdot
1/U^{1/(n-(n-k))}.$ It follows from Theorem 1 (letting $k_*=n-k$)
that $Y^*$ has a Pareto distribution $G_*(y)=\lambda_*
y/(1+\lambda_* y),\ y\ge 0,$ for some constant $\lambda_*>0.$
 This in turn implies that $Y$ has the Pareto
distribution $G(y)=\lambda y/(1+\lambda y),\ y\ge 0,$ where
$\lambda=1/\lambda_*.$ We can prove the last claim directly, or by
using Lemma 7(i), because $\lambda_* Y^*\stackrel{\rm
d}{=}1/(\lambda_* Y^*)=Y/\lambda_*$ having the standard Pareto type
II distribution. The proof is complete.
\medskip\\
\noindent{\bf Theorem 2.} Let $Y\sim G$ be a positive random
variable and let $1\le k \le n-1$ be fixed integers.
 Let $Y_1,Y_2,\ldots,Y_n$ be $n$ independent copies of $Y,$ and let $B$, independent of $\{Y_i\}_{i=1}^n$, be
a random variable having beta distribution $B_{\alpha,\beta}$  with
parameters $\alpha=1$ and $\beta=n-k,$ that is,
$F_B(u)=1-(1-u)^{n-k},\ u\in[0,1].$ Assume further that $\lim_{y\to
0^+}G(y)/y=\lambda>0.$ Then the distributional  equality
\begin{eqnarray}
Y_{k,k}\stackrel{\rm d}{=}Y_{k,n}/B
\end{eqnarray}
holds  iff $G$ is the Pareto  distribution
$G(y)=\lambda y/(1+\lambda y),\ y\ge 0.$\\
{\bf Proof.} The sufficiency part follows from Lemma 7(iii) with
$n-m=k$ and the fact $B\stackrel{\rm d}{=}U_{1,n-k}$. To prove the
necessity part, we note first that $G$ is absolutely continuous by
Lemmas 1--3, and then rewrite (7) as the functional equation:
\begin{eqnarray}
G^k(y)=\int_0^1\int_0^{G(yu)}k\bin{n}{k}t^{k-1}(1-t)^{n-k}dtdF_B(u),\ \ y\ge 0.
\end{eqnarray}
Now, it suffices to prove that the solution of equation (7) is
unique under the smoothness condition on the distribution. Namely,
if the absolutely continuous distribution $F$ on $(0,\infty)$
satisfies $\lim_{y\to 0^+}F(y)/y=\lambda>0$ and
\begin{eqnarray}
F^k(y)=\int_0^1\int_0^{F(yu)}k\bin{n}{k}t^{k-1}(1-t)^{n-k}dtdF_B(u),\ \ y\ge 0,
\end{eqnarray}
then we will prove that $F=G.$ From (8) and (9) it follows that
\begin{eqnarray}
|F^k(y)-G^k(y)|\le \frac{1}{E[B^k]}\int_0^1|F^k(yu)-G^k(yu)|dF_B(u),\ \ y\ge 0,
\end{eqnarray}
where $E[B^k]=1/\bin{n}{k}.$ Define the bounded function
$$g(y)=\left|\frac{F^k(y)-G^k(y)}{y^k}\right|,\ y>0,~~\hbox{and}~~g(0)=\lim_{y\to 0^+}g(y)=0,$$
and the increasing function
$$h(y)=\sup_{0\le t\le y}g(t),\ y>0,~~\hbox{and}~~h(0)=\lim_{y\to 0^+}h(y)=0.$$
By (10), we see that
\begin{eqnarray}
g(y)\le\int_0^1g(uy)dH(u),\ \ y\ge 0,
\end{eqnarray}
where $H(u)=(1/E[B^k])\int_0^ut^kdF_B(t),\ u\in[0,1].$ Now, by (11)
and the definition of the increasing function $h$, we have
\begin{eqnarray*}
h(y)\le\int_0^1h(uy)dH(u)\le h(y)\int_0^1dH(u)=h(y),\ \ y\ge 0.
\end{eqnarray*}
This in turn implies that $h$ is a constant function and hence
$h(y)=0,\ y\ge 0,$ because $h(0)=h(0^+)=0.$ Consequently, $g(y)=0,\
y\ge 0,$ and $F=G.$ The proof is complete.

The next result is the counterpart of Theorem 2 for the minimum
order statistics.

\noindent{\bf Corollary 3.} Under the same setting in Theorem 2 with
the condition on $G$ replaced by $\lim_{y\to
0^+}\overline{G}(1/y)/y=1/\lambda>0$ (equivalently,  $\lim_{y\to
+\infty}y\overline{G}(y)=1/\lambda>0$),
 the distributional equality
$Y_{1,k}\stackrel{\rm d}{=}Y_{n-k+1,n}\cdot B$ holds  iff $G$ is the
Pareto distribution $G(y)=\lambda y/(1+\lambda y),\ y\ge
0.$\\
{\bf Proof.} Use Lemma 7(v), Theorem 2 and the fact $(1/Y)_{\ell,
m}\stackrel{\rm d}{=}1/ Y_{m-\ell+1,m}$ for all $1\le\ell\le m.$

We now further extend Theorem 2 under some stronger smoothness
conditions.

\noindent{\bf Theorem 3.}  Let $Y\sim G$ be a positive random
variable and let $n,m,k$ be three fixed positive integers with $1\le
k\le n-m.$ Let $Y_1,Y_2,\ldots,Y_n$ be $n$ independent copies of
$Y,$
 and let $B_1$, independent of
$\{Y_i\}_{i=1}^n$, be a random variable having beta distribution
$B_{\alpha,\beta}$ with parameters $\alpha=n-m-k+1$ and $\beta=m.$
Assume further that the distribution function $G$
satisfies the following conditions:\\
(i) $G$ is real analytic and strictly increasing in $[0,\infty)$ and
 for each $i\geq 1$, its $i$-th derivative $G^{(i)}$ is strictly
monotone in some interval $[0,\delta_i)$.\\
(ii) $\lim_{y\to 0^+}[G^k(y)-(\lambda y)^k]/(\lambda y)^{k+1}=-k$
for some positive constant $\lambda.$\\ Then the distributional
equality
\begin{eqnarray}
Y_{k,n-m}\stackrel{\rm d}{=}Y_{k,n}/B_1
\end{eqnarray}
holds iff $G$ is the Pareto distribution
$G(y)=\lambda y/(1+\lambda y),\ y\ge 0.$\\
{\bf Proof.} The sufficiency part follows from Lemma 7(iii) and the
fact $B_1\stackrel{\rm d}{=} U_{n-m-k+1,n-k}.$ To prove the
necessity part, we note first that $G$ is absolutely continuous as
before, and then we rewrite (12) as the functional equation:
\begin{eqnarray}
G_{k,n-m}(y)=\int_0^1G_{k,n}(yu)dF_{B_1}(u),\ \ y\ge 0.
\end{eqnarray}
Now, it suffices to prove that the solution of equation (12) is
unique under the smoothness condition on the distribution. Namely,
if  the absolutely continuous distribution $F$ on $(0,\infty)$
satisfies the above conditions (i) and (ii) and
\begin{eqnarray}
F_{k,n-m}(y)=\int_0^1F_{k,n}(yu)dF_{B_1}(u),\ \ y\ge 0,
\end{eqnarray}
then we will prove that $F=G.$

Define the increasing function $H(y)=\max\{F(y), G(y)\},\ y\ge 0.$
From (1) it follows that for any $a>0$ and $0\le y\le a$,
\begin{eqnarray}
\left|F_{k,n-m}(y)-G_{k,n-m}(y)\right|&=&k\bin{n-m}{k}\left|\int_{G(y)}^{F(y)}t^{k-1}(1-t)^{n-m-k}dt\right|\nonumber\\
&\ge &k\bin{n-m}{k}(1-H(y))^{n-m-k}\frac{1}{k}\left|F^k(y)-G^k(y)\right|\nonumber\\
&\ge &\bin{n-m}{k}(1-H(a))^{n-m-k}\left|F^k(y)-G^k(y)\right|.
\end{eqnarray}
On the other hand, we have
\begin{eqnarray}
\left|F_{k,n}(y)-G_{k,n}(y)\right|\le \bin{n}{k}\left|F^k(y)-G^k(y)\right|,\ \ y\ge 0.
\end{eqnarray}
Combing (12)--(16) leads to
\begin{eqnarray}
& &\bin{n-m}{k}(1-H(a))^{n-m-k}\left|F^k(y)-G^k(y)\right|\le \left|F_{k,n-m}(y)-G_{k,n-m}(y)\right|\nonumber\\
&\le & \int_0^1 \left|F_{k,n}(yu)-G_{k,n}(yu)\right|dF_{B_1}(u)\le \bin{n}{k}\int_0^1\left|F^k(yu)-G^k(yu)\right|dF_{B_1}(u).
\end{eqnarray}

Define the bounded increasing function
$$h(y)=\sup_{0< t\le y}\left|\frac{F^k(t)-G^k(t)}{t^{k+1}}\right|,\ \ y>0,~~~ \hbox{and}~~h(0)=\lim_{y\to 0^+}h(y)=0.$$
Then from the inequality (17) it follows that for any $a>0$,
\begin{eqnarray}
& &\bin{n-m}{k}(1-H(a))^{n-m-k}h(y)\le \bin{n}{k}\int_0^1 h(yu)u^{k+1}dF_{B_1}(u)\nonumber\\
&\le & \bin{n}{k}h(y)\int_0^1u^{k+1}dF_{B_1}(u)=\bin{n}{k}h(y)E[B_1^{k+1}],\ \ 0\le y\le a.
\end{eqnarray}
Recall that $E[B_1^k]=\bin{n-m}{k}/\bin{n}{k}.$ Then rewrite the
inequality (18) as follows:
\begin{eqnarray}
& &E[B_1^k](1-H(a))^{n-m-k}h(y)\le h(y)E[B_1^{k+1}],\ \ 0\le y\le a,\ a>0.
\end{eqnarray}
We claim that there exists a $y_0>0$ such that $h(y_0)=0$.
Otherwise, we have, by (19),
\begin{eqnarray*}
& &E[B_1^k](1-H(a))^{n-m-k}\le E[B_1^{k+1}],\ \ \forall\ a>0,
\end{eqnarray*}
which in turn implies, by letting $a\to 0^+$, that $E[B_1^k]\le
E[B_1^{k+1}]$, a contradiction. Therefore, $h(y_0)=0$ for some
$y_0>0$ and hence $F(y)=G(y)$ for $y\in [0, y_0]$. By Lemma 9 and
the assumptions on $F$ and $G$, we conclude that $F=G$. The proof is
complete.

\noindent{\bf Corollary 4.}   Let $Y\sim G$ be a positive random
variable and  let $n,m,k$ be three fixed positive integers with
$k+1\le m\le n.$ Let $Y_1,Y_2,\ldots,Y_n$ be $n$ independent copies
of $Y,$
 and  let $B_2$, independent of
$\{Y_i\}_{i=1}^n$, be a random variable having beta distribution
$B_{\alpha,\beta}$ with parameters $\alpha=m-k$ and $\beta=k$.
 Assume
further that the distribution function $G_*$ of $1/Y$
satisfies the following conditions:\\
(i) $G_*$ is real analytic and strictly increasing in $[0,\infty)$
and
 for each $i\geq 1$, its $i$-th derivative $G_*^{(i)}$ is strictly
monotone in some interval $[0,\delta_i)$.\\
(ii) $\lim_{y\to
0^+}[G_*^{k_*}(y)-(y/\lambda)^{k_*}]/(y/\lambda)^{k_*+1}=-k_*$
for some positive constant $\lambda,$ where $k_*=n-m+1.$\\
 Then the distributional equality
$Y_{m-k,n-k}\stackrel{\rm d}{=}Y_{m,n}\cdot B_2$ holds iff $G$ is
the Pareto distribution
$G(y)=\lambda y/(1+\lambda y),\ y\ge 0.$ \\
{\bf Proof.} Use Lemma 7(v), Theorem 3  and the fact $(1/Y)_{\ell,
m}\stackrel{\rm d}{=}1/ Y_{m-\ell+1,m}$ for all $1\le\ell\le m.$

\indent In summary, for a positive random variable $Y\sim G,$ we
have the following characteristic properties of the Pareto
distribution $G(y)=\lambda y/(1+\lambda y),\ y\ge 0,$ where
$\lambda$ is a positive constant (compare with Lemma 7).
\begin{enumerate}
\item $Y_{k,n-1}\stackrel{\rm d}{=}Y_{k,n}/U^{1/(n-k)}.$
\item $Y_{k,n-1}\stackrel{\rm d}{=}Y_{k+1,n}\cdot U^{1/k}.$
\item $Y_{k,k}\stackrel{\rm d}{=}Y_{k,n}/B$ (equivalently, $Y_{k,k}\stackrel{\rm d}{=}Y_{k,n}/U_{1, n-k}$).
\item $Y_{1,k}\stackrel{\rm d}{=}Y_{n-k+1,n}\cdot B$ (equivalently, $Y_{1,k}\stackrel{\rm d}{=}Y_{n-k+1,n}\cdot U_{1, n-k}$).
\item $Y_{k,n-m}\stackrel{\rm d}{=}Y_{k,n}/B_1$ (equivalently, $Y_{k,n-m}\stackrel{\rm d}{=}Y_{k,n}/U_{n-m-k+1,n-k}$).
\item $Y_{m-k,n-k}\stackrel{\rm d}{=}Y_{m,n}\cdot B_2$ (equivalently, $Y_{m-k,n-k}\stackrel{\rm d}{=}Y_{m,n}\cdot  U_{m-k,m-1}$).
\end{enumerate}
Here, the random variables $U\sim {\cal U},\ B\sim B_{1, n-k},\
B_1\sim B_{n-m-k+1, m}$, $B_2\sim B_{m-k, k},$ and on the RHS of
each equality, the two random variables are independent.
\medskip\\
 \noindent{\bf 4. Characterizations of the logistic
distribution}

\indent We are now ready to provide characterization results of the
logistic distribution.

\noindent{\bf Theorem 4.} Let $X\sim F$ and let $1\le k \le n-1$ be
fixed integers. Then the distributional  equality
\begin{eqnarray*}
X_{k,n-1}\stackrel{\rm d}{=}X_{k,n}+\frac{1}{n-k}\xi
\end{eqnarray*}
holds  iff $F$ is a logistic distribution
$F(x)=1/[1+e^{-(x-\mu)}],\ x\in{\mathbb R},$  where $\mu\in{\mathbb R}$ is a constant.\\
{\bf Proof.} Take $Y_i=\exp(X_i)$, $U=\exp(-\xi)$ and
$\lambda=e^{-\mu}.$ Then the result follows immediately from Theorem
1.

\noindent{\bf Corollary 5.} Let $X\sim F,$  $\alpha>0$ and let $1\le
k \le n-1$ be fixed integers.   Then the distributional equality
$X_{k,n-1}\stackrel{\rm d}{=}X_{k,n}+\frac{1}{\alpha}\xi$ holds iff
$F$ is a logistic distribution
$F(x)=1/\{1+e^{-[\alpha/(n-k)](x-\mu)}\},\ x\in{\mathbb R},$ where
$\mu\in{\mathbb R}$ is a constant.

\noindent{\bf Corollary  6.}  Let $X\sim F,$  $\alpha>0$ and  let
$1\le k \le n-1$ be fixed integers. Then  the distributional
equality
\begin{eqnarray}
X_{k,n-1}\stackrel{\rm d}{=}X_{k+1,n}-\frac{1}{\alpha}\xi
\end{eqnarray}
holds  iff  $F$ is a logistic distribution
$F(x)=1/[1+e^{-(\alpha/k)(x-\mu)}],\ x\in{\mathbb R},$
where $\mu\in{\mathbb R}$ is a constant.\\
{\bf Proof.} Use Corollary 5 and the fact that
$X_{\ell,m}\stackrel{\rm d}{=}-(-X)_{m-\ell+1,m}$ for all $1\le
\ell\le m.$

The counterpart of (20), namely, $X_{k+1,n}\stackrel{\rm
d}{=}X_{k,n-1}+\frac{1}{\alpha}\xi,$ and the two-sided case:
$X_{k,n}+a\xi_1\stackrel{\rm d}{=}X_{k,n-1}+b\xi_2$ (see Corollary
5), where $\alpha, a, b>0$, were investigated by
Weso$\hbox{\l}$owski and Ahsanullah (2004). All the solutions of
these two equations are exponential distributions.

The next result improves and extends Theorem 6 of Lin and Hu (2008)
by an approach different from the previous method of intensively
monotone operator (Kakosyan et al.\,1984).

\noindent{\bf Theorem 5.} Let  $1\le k \le n-1$ be fixed integers.
Assume that $X\sim F$ satisfies $\lim_{x\to
-\infty}F(x)/e^x=e^{-\mu}$ for some constant $\mu\in{\mathbb R}.$
 Then the distributional   equality
\begin{eqnarray*}
X_{k,k}\stackrel{\rm d}{=}X_{k,n}+\xi_{n-k,n-k}
\end{eqnarray*}
holds iff $F$ is the logistic distribution
$F(x)=1/[1+e^{-(x-\mu)}],\ x\in{\mathbb R}.$\\
{\bf Proof.} The sufficiency part follows from Lemma 8(ii), while
the necessity part is a consequence of Theorem 2.

The next result is the counterpart of Theorem 5 for the minimum
order statistics.

\noindent{\bf Corollary 7.} Let $1\le k \le n-1$ be fixed integers.
 Assume that $X\sim F$
satisfies $\lim_{x\to -\infty}\overline{F}(-x)/e^x=e^{\mu}$ for some
constant $\mu\in{\mathbb R}.$ Then the distributional equality
\begin{eqnarray*}
X_{1,k}\stackrel{\rm d}{=}X_{n-k+1,n}-\xi_{n-k,n-k}
\end{eqnarray*}
holds iff $F$ is the logistic distribution
$F(x)=1/[1+e^{-(x-\mu)}],\ x\in{\mathbb
R}.$\\
{\bf Proof.} Use Lemma 8(iv), Theorem 5 and the fact
 $X_{\ell,m}\stackrel{\rm
d}{=}-(-X)_{m-\ell+1,m}$ for all $1\le \ell\le m.$

 Using Theorem 3 and Lemma 8(ii), we further extend Theorem 5 to the following.

\noindent{\bf Theorem 6.} Let $n,m,k$ be three fixed positive
integers with $1\le k\le n-m$ and let $X\sim F$ satisfy $\lim_{x\to
-\infty}[e^{-k(x-\mu)}F^k(x)-1]/e^{x-\mu}=-k$  for some constant
$\mu\in{\mathbb R}.$ Assume further that the distribution function
$G$ of $\exp(X_1)$ is real analytic and strictly increasing in
$[0,\infty)$ and that
 for each $i\geq 1$, its $i$-th derivative $G^{(i)}$ is strictly
monotone in some interval $[0,\delta_i)$. Then the distributional
equality
\begin{eqnarray*}
X_{k,n-m}\stackrel{\rm d}{=}X_{k,n}+\xi_{m,n-k}
\end{eqnarray*}
holds iff $F$ is the logistic distribution
$F(x)=1/[1+e^{-(x-\mu)}],\ x\in{\mathbb R}.$

As before, Theorem 6 and Lemma 8(iv) together lead to the following.

\noindent{\bf Corollary 8.}  Let $n,m,k$ be three fixed positive
integers with $k+1\le m\le n$ and let $X\sim F$ satisfy $\lim_{x\to
-\infty}[e^{-k_*(x+\mu)}(\overline{F}(-x))^{k_*}-1]/e^{x+\mu}=-k_*$
for some constant $\mu\in{\mathbb R},$ where $k_*=n-m+1.$ Assume
further that the distribution function $G_*$ of $\exp(-X_1)$ is real
analytic and strictly increasing in $[0,\infty)$ and that
 for each $i\geq 1$, its $i$-th derivative $G_*^{(i)}$ is strictly
monotone in some interval $[0,\delta_i)$. Then the distributional
equality
\begin{eqnarray}
X_{m-k,n-k}\stackrel{\rm d}{=}X_{m,n}-\xi_{k,m-1}
\end{eqnarray}
holds iff $F$ is the logistic distribution
$F(x)=1/[1+e^{-(x-\mu)}],\ x\in{\mathbb R}.$

\indent In summary, for a random variable $X\sim F,$ we have the
following characteristic properties of the logistic distribution
$F(x)=1/[1+e^{-(x-\mu)}],\ x\in{\mathbb R}$  (compare with Lemma 8).
Here, $\mu\in{\mathbb R},$ $\xi\sim {\cal E},$ and on the RHS of
each equality, the two random variables are independent.
\begin{enumerate}
\item $X_{k,n-1}\stackrel{\rm d}{=}X_{k,n}+\frac{1}{n-k}\xi.$
\item $X_{k,n-1}\stackrel{\rm d}{=}X_{k+1,n}-\frac{1}{k}\xi.$
\item $X_{k,k}\stackrel{\rm d}{=}X_{k,n}+\xi_{n-k,n-k}.$
\item $X_{1,k}\stackrel{\rm d}{=}X_{n-k+1,n}-\xi_{n-k,n-k}.$
\item $X_{k,n-m}\stackrel{\rm d}{=}X_{k,n}+\xi_{m,n-k}.$
\item $X_{m-k,n-k}\stackrel{\rm d}{=}X_{m,n}-\xi_{k,m-1}.$
\end{enumerate}
\noindent{\bf 5. Open problem}\\
\indent Finally, we would like to pose the following problem in
which part (i) is the counterpart of (21) for exponential
distribution, and  in part (ii), the first two cases, $m=k+1,k+2$,
have been solved by AlZaid and Ahsanullah (2003) and Ahsanullah et
al.\,(2010).
\medskip\\
\noindent {\bf  Problem.} Let $X\sim F$ and let $1\le k<m\le n$ be
fixed integers. Then solve the general distributional equations: (i)
$X_{m,n}\stackrel{\rm d}{=}X_{m-k,n-k}+\xi_{k,n}$ and (ii)
$X_{m,n}\stackrel{\rm d}{=}X_{k,n}+\xi_{m-k,n-k}.$
\medskip\\
\noindent{\bf Acknowledgments.} The authors would like to thank the
Editor-in-Chief and the Referee for helpful comments and
constructive suggestions which improve the presentation of the
paper.\bigskip\\
 \centerline{\bf References}
\begin{description}

\item Ahsanullah, M.,  Berred, A. and Nevzorov, V.B. (2011).  On characterizations
of the exponential distributions.  {\it J. Appl. Statist. Sci.},
{\bf 19}, 37--43.

\item Ahsanullah, M., Nevzorov, V.B. and Yanev, G.P. (2010). Characterizations of distributions
via order statistics with random exponential shifts. {\it J.
Appl. Statist. Sci.}, {\bf 18}, 297--305.

\item Ahsanullah, M., Yanev, G.P. and Onica, C. (2012). Characterizations of logistic distribution
through order statistics with independent exponential shifts.
{\it Economic Quality Control}, {\bf 27}, 85--96.

\item AlZaid, A.A. and Ahsanullah, M. (2003). A characterization of the Gumbel distribution based on record values.
{\it Commun. Statist. - Theory and Methods}, {\bf 32},
2101--2108.

\item Ananjevskii, S.M. and Nevzorov, V.B. (2016). On families of
distributions characterized by certain properties of ordered
random variables. {\it Vestnik St. Petersburg University:
Mathematics}, {\bf 49}, 197--203.

\item Berred, A. and Nevzorov, V.B. (2013). Characterizations of families of distributions, which
include the logistic one, by properties of order statistics.
{\it Journal of Mathematical Sciences}, {\bf 188}, 673--676.

\item David, H.A. and Nagaraja, H.N. (2003). {\it Order Statistics}, 3rd edn.
Wiley, New Jersey.

\item David, H.A. and Shu, V.S. (1978). Robustness of location estimators
in presence of an outlier. In: {\it Contributions to Survey
Sampling and Applied Statistics}, ed. H.A. David, pp.\,235--250.
Academic Press, New York.

\item George, E.O. and Mudholkar, G.S. (1982). On the logistic
and exponential laws. {\it Sankhy$\bar{a}$, Ser. A}, {\bf 44}, 291--293.

\item Gradshteyn, I.S. and Ryzhik, I.M. (2014). {\it Table of Integrals, Series,
and Products},  8th edn.  Elsevier, New York.

\item Hwang, J.S. and Lin, G.D. (1984). Characterizations of
distributions by linear combinations of moments of order
statistics. {\it Bulletin of the Institute of Mathematics,
Academia Sinica}, {\bf 12}, 179--202.

\item Kakosyan, A.V., Klebanov, L.B. and Melamed, J.A. (1984).
{\it Characterization of Distributions by the Method of
Intensively Monotone Operators.} Springer, New York.

\item Kotz, S. (1974). Characterizations of statistical distributions:
a supplement to recent surveys. {\it Int. Statist. Rev.},
{\bf 42}, 39--65.

\item Lin, G.D. and Hu, C.-Y. (2008). On characterizations of the logistic
distribution. {\it J. Statist. Plann. Infer.}, {\bf 138},
1147--1156.

\item Lukacs, E. (1970). {\it Characteristic Functions}, 2nd edn. Hafner Pub. Co., New York.

\item Mudholkar, G.S. and George, E.O. (1978). A remark on the shape
of the logistic distribution. {\it Biometrika}, {\bf 65}, 667--668.

\item Nevzorov, V.B. (2001). {\it Records: Mathematical Theory}. Translations
of Mathematical Monographs, Vol. 194, Amer. Math. Soc., Rhode Island.

\item Shorack, G.R. and Wellner, J.A. (1986). {\it Empirical Processes with Applications to Statistics.} Wiley, New York.

\item Weso$\hbox{\l}$owski, J. and Ahsanullah, M. (2004). Switching order statistics through random power contractions.
{\it Aust. N.Z. J. Statist.}, {\bf 46}, 297--303.

\item Zykov, V.O. and Nevzorov, V.B. (2011). Some characterizations of
families of distributions, including logistic or exponential
ones, by properties of order statistics. {\it Journal of
Mathematical Sciences}, {\bf 176}, 203--206.
\end{description}
\end{document}